\documentclass[11pt,a4paper]{article}
\usepackage{graphicx}
\usepackage{latexsym}
\usepackage[all]{xy}
\usepackage{amsfonts}
\usepackage{amsthm}
\usepackage{amsmath}
\usepackage{amssymb}

\def\<{{\langle}}
\def\>{{\rangle}}

\def\note#1{{}}

\def\note#1{}

\def\beq{\begin{equation}}
\def\eeq{\end{equation}}

\def\ot{{\otimes}}

\newcounter{zlist}

\def\Label#1{\label{#1}\ifmmode\llap{[#1] }\else 
\marginpar{\smash{\hbox{\tiny [#1]}}}\fi}
\def\Label{\label}

\newtheorem{proposition}{Proposition}[section]

\newtheorem{theorem}[proposition]{Theorem}

\theoremstyle{definition}
\newtheorem{definition}[proposition]{Definition}

\theoremstyle{remark}
\newtheorem{remark}[proposition]{Remark}

\newcounter{c}

\newcommand{\etyk}[1]{\vspace{-7.4mm}$$\begin{equation}\Label{#1}
\addtocounter{c}{1}}
\renewcommand{\]}{\ifnum \value{c}=1 $$\else \end{equation}\fi}
\begin{document}

\title{ \bf NON-ASSOCIATIVE ALGEBRAS,
YANG-BAXTER EQUATIONS  AND QUANTUM COMPUTERS}

\author{RADU IORDANESCU\\
Institute of Mathematics "Simion Stoilow" of the Romanian Academy \\ 
P.O. Box 1-764, RO-014700 Bucharest, Romania \\
{ \it E-mail: Radu.Iodanescu@imar.ro}\\
FLORIN F. NICHITA\\
Institute of Mathematics "Simion Stoilow" of the Romanian Academy \\ 
P.O. Box 1-764, RO-014700 Bucharest, Romania \\
{ \it E-mail: Florin.Nichita@imar.ro}\\
ION M. NICHITA\\
ASESOFT Ploiesti, Romania\\
{ \it E-mail: mmiihhuu@yahoo.com}
}

\date{}
\maketitle
\begin{abstract}
Non-associtive algebras is a research direction 
gaining much attention these days.
New developments show that associative algebras and some not-associative
 structures can be unified at the level of Yang-Baxter structures.
In this paper, we present a unification for associative algebras,
Jordan algebras and Lie algebras.
 
The (quantum) Yang-Baxter equation and related structures
are interesting topics, because they
have applications in many
areas of mathematics,
physics and computer science. 
Several new interpretations and results are presented below.

\end{abstract}

\bigskip

{\bf AMS 2010 Subj. Class.:}  17C05, 17C50, 16T15, 16T25, 17B01

\bigskip


\bigskip

{\centering\section{INTRODUCTION}}

The current paper emerged after our participation at the 
International Conference
"Mathematics Days in Sofia", July 7-10, 2014, Sofia, Bulgaria.
It  presents our scientic contributions for that conference,
as well as new results and directions of research.

Prof. Radu Iordanescu contributed with a paper on Jordan algebras (\cite{I}) to the library of the Institute of Mathematics and Informatics from Sofia.
Non-associtive algebras (which include Jordan algebras, Lie algebras and associative algebras)
are currently a research direction in fashion (see \cite{wt}, and the references therein).
New developments show that associative algebras and Lie algebras can be unified at the level of Yang-Baxter structures (see \cite{humbold}).
In this paper, we present a unification for associative algebras,
Jordan algebras and Lie algebras.

The other authors of this paper presented a poster on combinatorial logical circuits
and solutions to the set-theoretical Yang-Baxter equation from Boolean algebras (following the works \cite{fm, f}).
Dr. Violeta Ivanova (see \cite{vio}) was interested in the applications of these problems in computer science.
The Yang-Baxter equation can be interpretated in terms of combinatorial logical structures, and, in logic, it 
represents some kind of compatibility condition, when working with many logical sentences in the same time.
This equation is also related to the theory of universal quantum gates and to the quantum computers (see, for example, \cite{abj}).

Florin F. Nichita gave a talk on the Yang-Baxter equation presenting resuts from \cite{ax, ax2}. In his talk he 
referred to several papers of Prof. Vladimir Dobrev and to the work of Prof. Tatiana Gateva-Ivanova (one of her
 questions, arising at that time,
will be partially answered in this paper). As an observation, our Yang-Baxter operator $ R^{A}_{\alpha, \beta, \alpha} $
is related to a universal quantum gate.

The organization of our paper is the following. The next section introduces the mathematical terminology needed 
in this paper,
and it presents results about the Yang-Baxter equation. Section 3 deals with the unification of
Jordan algebras, Lie algebras and associative algebras. 
Section 4 is a conlusions section and an update on quantum computers.


\bigskip

{\centering\section{THE QYBE}}

An introduction to the (quantum) Yang-Baxter equation (QYBE) 
could be found in the paper \cite{ax}.
Several special sessions on it followed at the open-acces journal AXIOMS,
explaining its role in
areas of mathematics,
physics and computer science. 

\bigskip

We will work over the field $ k $, and the tensor products will be defined over $k$.
For $ V $ a $ k$-space, we denote by
$ \   \tau : V \otimes V \rightarrow V \ot V \  $ the twist map defined by $ \tau (v \ot w) = w \ot v $, and by $ I: V \rightarrow V $
the identity map of the space V.
For $ \  R: V \ot V \rightarrow V \ot V  $
a $ k$-linear map, let
$ {R^{12}}= R \ot I , \  {R^{23}}= I \ot R , \
{R^{13}}=(I\ot\tau )(R\ot I)(I\ot \tau ) $.

\bigskip

\begin{definition} A
{ \it  Yang-Baxter
operator} 
is defined as an invertible  $ k$-linear map  $ R : V \ot V \rightarrow V \ot V $
which satisfies the  equation:
\begin{equation}  \label{ybeq}
R^{12}  \circ  R^{23}  \circ  R^{12} = R^{23}  \circ  R^{12}  \circ  R^{23}
\end{equation}
$R$ satisfies (\ref{ybeq}) if and only if
$R\circ \tau  $ ( respectively $ \tau \circ R $) satisfies the QYBE:
\begin{equation}   \label{ybeq2}
R^{12}  \circ  R^{13}  \circ  R^{23} = R^{23}  \circ  R^{13}  \circ  R^{12}
\end{equation}
\end{definition}

\bigskip
\begin{remark} \label{alg}
For $A$ be a (unitary) associative $k$-algebra, and $ \alpha, \beta, \gamma \in k$,
  \cite{DasNic:yan} defined the
$k$-linear map
$$   R^{A}_{\alpha, \beta, \gamma}: A \ot A \rightarrow A \ot A, \ \ 
R^{A}_{\alpha, \beta, \gamma}( a \ot b) = \alpha ab \ot 1 + \beta 1 \ot ab -
\gamma a \ot b \ , $$\\
which is a Yang-Baxter operator if and only if one
of the following holds:\\ 
(i) $ \alpha = \gamma \ne 0, \ \ \beta \ne 0 $; $ \ $
(ii) $ \beta = \gamma \ne 0, \ \ \alpha \ne 0 $; $ \ $
(iii) $ \alpha = \beta = 0, \ \ \gamma \ne 0 $.
\end{remark}

\begin{remark}
Using
Turaev's general scheme to derive an invariant of 
oriented links from a Yang-Baxter operator  (see \cite{T}),
$ \  R^{A}_{\alpha, \beta, \gamma} $ leads to the 
Alexander polynomial of knots (see \cite{mn}).

\end{remark}

\bigskip

\begin{remark}
In dimension two, $ R^{A}_{\alpha, \beta, \alpha} \circ \tau$,
can be expressed as:
\begin{equation} \label{rmatcon2}
\begin{pmatrix}
1 & 0 & 0 & 0\\
0 & 1 & 0 & 0\\
0 & 1-q  & q & 0\\
\eta & 0 & 0 & -q
\end{pmatrix}
\end{equation}
where $ \eta \in \{ 0, \ 1 \} $, and $q \in k - \{ 0 \}$. For  $ \eta = 0 $ and $q = 1$,
it can be related to the universal quantum gate CNOT.
\end{remark}

\bigskip

\begin{definition}
A Lie superalgebra is a (nonassociative) $Z_2$-graded algebra, or superalgebra, 
over a field $k$ with the  Lie superbracket, satisfying the two conditions:
$$[x,y] = -(-1)^{|x||y|}[y,x$$
$$ (-1)^{|z||x|}[x,[y,z]]+(-1)^{|x||y|}[y,[z,x]]+(-1)^{|y||z|}[z,[x,y]]=0 $$
where $x$, $y$ and $z$ are pure in the $Z_2$-grading. Here, $|x|$ denotes the degree of $x$ (either 0 or 1). 
The degree of $[x,y]$ is the sum of degree of $x$ and $y$ modulo $2$.
\end{definition}

\bigskip

\begin{remark} \label{Lie}
For $ ( L , [,] )$  a Lie superalgebra over $k$,
  $ z \in Z(L) = \{ z \in L : [z,x]=0 \ \ \forall \ x \in L \} ,
 \ \vert z \vert =0 $ and $ \alpha \in k $, \cite{mj} (and \cite{nipo})  defined
$$ \ \ \  { \phi }^L_{ \alpha} \ : \ L \ot L \ \ \longrightarrow \ \  L \ot L, \ \ 
x \ot y \mapsto \alpha [x,y] \ot z + (-1)^{ \vert x \vert \vert y \vert } y \ot x \ , $$
which is a YB operator.
\end{remark}

\bigskip

\begin{remark} \label{unif}
 The Remarks \ref{alg} and \ref{Lie} lead to some kind of unification of
associative algebras and structures that are not associative at the level
of Yang-Baxter structures (see \cite{humbold,nichita}. For example, the first
isomorfism theorem for groups (algebras) and the first
isomorfism theorem for
Lie algebras, can be unified as
an isomorphism theorem for Yang-Baxter structures (see \cite{nonlin}).
The fact that Therem 7.2.3 from \cite{nonlin} (The fundamental
isomorphism theorem for YB structures)
 unifies not only associative algebras and coalgebras, but also
Lie algebras is a new result.

\end{remark}

\bigskip

\begin{remark} Following a question of  Prof. Tatiana Gateva-Ivanova,
 we can construct an algebra structure associated to the operator\\
$   R= R^{A}_{1, 1, 1}: A \ot A \rightarrow A \ot A, \ \ 
R( a \ot b) =  ab \ot 1 +  1 \ot ab -
a \ot b \ , $\\
if we use Theorem 3.1 (i) from \cite{ax2}.


For $ a, b \in T^1 (A)=A $, we have:
$ \mu(a \ot b) = ab \ot 1 + 1 \ot ab - a \ot b \ \in T^2 (A) \ . $


 For $ a \ot a' \in T^2 (A)= A \ot A $ and  $ b \in T^1 (A)=A $, we have:

$ \mu((a \ot a') \ot b) = R^{12} \circ R^{23} (a \ot a' \ot b) =
aa'b \ot 1 \ot 1 + 1 \ot aa'b \ot 1 - a \ot a'b \ot 1
+ 1 \ot a \ot a'b -
aa' \ot 1 \ot b - 1 \ot aa' \ot b + a \ot a' \ot b 
\in T^3 (A) \ . $


 For $ a \in T^1 (A)= A  $ and  $ b \ot b' \in T^2 (A) $, we have:

$ \mu(a \ot ( b \ot b') = R^{23} \circ R^{12} (a \ot b \ot b') \in T^3 (A) \ . $

In the same manner, we compute other products.

\end{remark}

\bigskip

{\centering\section{NON-ASSOCIATIVE ALGEBRAS}}

Jordan algebras emerged in the early thirties, and their applications
are in differential geometry, ring geometries, physics, quantum groups,
analysis, biology, etc (see \cite{I, RI}). A poster
presented at the 11-th International Workshop on Differential Geometry and its 
Applications, in September 2013, at the Petroleum-Gas University
from Ploiesti led to the paper \cite{florin}.
That paper presents new results on Jordan algebras, Jordan 
coalgebras, and how they are related
to the QYBE. 

\bigskip

One of the main results of \cite{florin} is the following theorem, which
explaines when the Jordan identity implies associativity. It is an intrinsic result.

\bigskip

\begin{theorem}
Let  V   be a vector space spanned by   $a$ and  $b$, which are linearly independent.
Let                                                                                 $ \theta : V \otimes V \rightarrow V, \ \ 
\theta (x \otimes y) = xy $, 
be a linear map which is a commutative operation
with the property

 \begin{equation} \label{sase}
a^2 = b \ , \ \ \ b^2 = a \ .
\end{equation} 
                                                                               
Then:
$ ( V, \ \theta ) $
     is a Jordan algebra  
$ \iff $
$ ( V, \ \theta ) $
      is a non-unital commutative (associative) algebra.
\end{theorem}

\bigskip

The next remark finds a relationship between Jordan algebras, Lie algebras and associative algebras.
In this case, we have an extrinsic result about non-associative structures.

\bigskip

\begin{remark}
For the vector space V, let 
   $ \eta : V \otimes V \rightarrow V, \ \ 
\eta (x \otimes y) = xy \ , \ \ $  be a linear map which satisfies:
 
\begin{equation} \label{new}
 (ab)c + (bc)a + (ca)b \ = \ a(bc) + b(ca) + c(ab) \ \ ;
\end{equation} 

\begin{equation} \label{Jordan}
 (a^2 b) a \ = \ a^2 (ba) \ \ .
\end{equation}

Then, $(V, \eta) $ is a structure which unifies (non-unital) associative algebras, Lie algebras
and Jordan algebras.

Indeed, the associativity and the Lie identity are unified by relation (\ref{new}). Also, the commutativity
of a Jordan algebra implies (\ref{new}).
But, the Jordan identity, (\ref{Jordan}), which appears in the definition of Jordan algebras,
holds true in any associative algebra and Lie algebra. 

This unifying approach is different from the unification proposed in Remark \ref{unif}.

\end{remark}

\bigskip

{\centering\section{CONCLUSIONS AND QUANTUM COMPUTERS}}

The first quantum computer (which uses principles of quantum mechanics) was sold to the aerospace and security of defense company Lockheed Martin.
The manufacturing company, D-Wave, founded in 1999 and called ``a company of quantum computing" promised to perform professional services for
the computer maintenance as well.

The Yang-Baxter equation has applications in quantum computing, and it can be viewed as an unifying structure for non-associative structures.

The authors of this paper would like to thank the organizers of the
International Conference
"Mathematics Days in Sofia", 2014, 
for their kind hospitality.




\begin{center}

\end{center}

\bigskip

\bigskip

\bigskip

\end{document}